\documentclass[11pt]{svjour3}          
\smartqed  
\usepackage{graphicx,epsfig}
\usepackage{amsmath}
\begin{document}

\title{Optimal Design of Helical Springs of Power Law Materials}

\author{Dongming Wei \and Marios Fyrillas \and \\ Adilet Otemissov \and Rustam Bekishev \\\\
Nazarbayev University,\\
 53 Kabanbay Batyr Ave., Astana, 010000, Kazakhstan}


\institute{Corresponding Author, \mbox{E-mail: dongming.wei@nu.edu.kz}\\
\mbox{E-mail: marios.fyrillas@nu.edu.kz}\\
\mbox{E-mail: aotemissov@nu.edu.kz}\\
\mbox{E-mail: Rustam.Bekishev@nu.edu.kz}}

\maketitle

\begin{abstract}
In this paper the geometric dimensions of a compressive helical spring made of power law materials are optimized to reduce the amount of material. The mechanical constraints are derived to form the geometric programming problem. Both the prime
and the dual problem are examined and solved semi-analytically for a range of spring index. A numerical example is provided to validate the solutions.
\end{abstract}

\begin{table}[h]
\begin{tabular}{|l|}  \hline
{\bf Nomenclature}  \\
$C$=spring index $=\frac{D}{d}$ \\
$d$=spring wire  diameter, $m$ \\
$D$=mean coil  diameter, $m$ \\
$N$=number of turns in the spring \\
$\delta$= tip deflection of the spring, $m$ \\
$K$=bulk modulus,  $MPa$ \\
$n$=the power law index \\
$\rho$=the density of the material, $kg/m^3$ \\
\hline
\end{tabular}
\end{table}

\section{Introduction}\label{sec:1}
\label{intro} Helical springs are the  basic structure elements used in many mechanical devices. In many applications, it is important to optimize the geometric dimensions of the springs to reduce the amount of material used while maintaining the ability to support the required loads. Optimal design of  helical springs based on geometric programing is well-known for materials which obey Hooke's law, see, e.g., \cite{AG}, \cite{SK}, \cite{MP}. However, materials subject to nonlinear stress-strain constitutive laws in this context have not been well-studied. One of the simplest nonlinear material constitutive law is the following power law \begin{equation}\label{powerlaw} \sigma=K|\epsilon|^{n-1}\epsilon\end{equation}
where $\sigma$ is the axial stress, $\epsilon$ the axial strain, $K$ the material constant-called the bulk modulus, and $n$ the power law index. Materials which obey (\ref{powerlaw}) are often called Hollomon or Ludwick materials in literature, see, e.g. \cite{WSE} and \cite{WL}.  High strength alloy metals such as heat treated metals, stainless steels, Titanium alloys, and the super plastic-polyimide are the common examples of power law materials.  See, e.g., \cite{KSK} and \cite{SC} for a list of common metals with numerical values of the bulk modulus and the power law index.

Helical springs are mechanical devices made of a wire coiled in the form of a helix  and considered to be a major element of shock absorber, return mechanisms, fuel flow controller  used in engineering, automotive, medical and agricultural machinery. They are widely used for compressive loads. Because even for small strains there is no obvious yield of the stress-strain curve for the power law materials before ultimate yield point, the linear theory or the traditional reduced modulus theory is not applicable or cannot be accurately applied to calculate stress distributions.

In this work, we first derive  the maximum mechanical stress or loads and the corresponding tip-deflection under a compressive load, and then formulate the corresponding geometric programming problem minimizing the amount of material needed for the given loads.  The optimal solutions of the geometric programming problem is studied by considering KKTC conditions. The corresponding dual problem for the primal problem is constructed and examined for the solution as well. A numerical example is provided for both the prime and the dual problem.

\section{The Mechanical Constraints for the Spring}\label{sec:2}
Let $x$ denote the distance along a circular shaft  from the fixed end under a uniform torque. We assume that the rotation at $x$, denoted by $\phi(x)$,  is proportional to $x$, i.e., $\phi(x)=\alpha x$, where $\alpha$ is the rate of twist. Further, we assume that $\epsilon_{xx}=\epsilon_{rr}=\epsilon_{xr}=\epsilon_{r\theta}=0$, and $\epsilon_{x\theta}=\frac{r\alpha}{2}$, where $\epsilon_{ij},i,j=x,r,\theta$ are the strains in polar coordinates. For the power law materials, we have the shear strain due to torsion $\tau_{x\theta}=2G|\epsilon_{x\theta}|^{n-1}\epsilon_{x\theta}=2G(\frac{r\alpha}{2})^n$, where $G=\frac{K}{1+\nu}$ is the shear modulus and $\nu$  the Poison's ratio.  Let $ A$ denote the cross-section of the shaft, the total torque at $x$ is given by $T=\int_A\tau_{x\theta}rdA=2G(\frac{\alpha}{2})^nI_n$, where $I_n=\int_A r^{n+1}drd\theta=\frac{\pi d^{n+3}}{(n+3)2^{n+2}}$ is the generalized area moment. Therefore, we have $\alpha =2(\frac{T}{2GI_n})^{1/n}$, $\tau_{x\theta}=\frac{Tr^n}{I_n}$. Assume that the angle between loading force at the tip of the spring  and the plane containing the cross-section $A$ of the wire is negligible, then the torque acting on the wire $T=P\frac{D}{2}$ and  the average shear stress in the wire due to a vertical constant load, denoted by $P$ is $\tau_{av} = \frac{Force}{Area} = \frac{P}{\pi d^2/4} = \frac{4P}{\pi d^2}$. Therefore, the total stress at $r$ is approximately
\begin{equation}\label{Stress}
\tau= \frac{Tr^n}{I_n}+\frac{4P}{\pi d^2}=\frac{2^{n+1}(n+3)PDr^n}{\pi d^{3+n}}+\frac{4P}{\pi d^2}
\end{equation}
We now derive the deflection of the tip by extending the standard textbook method ( see, e.g., \cite{Sp} and \cite{Sh} )for $n=1$ to the case of any value of $n$.
\begin{figure}[ht]
\centerline{\epsfig{figure=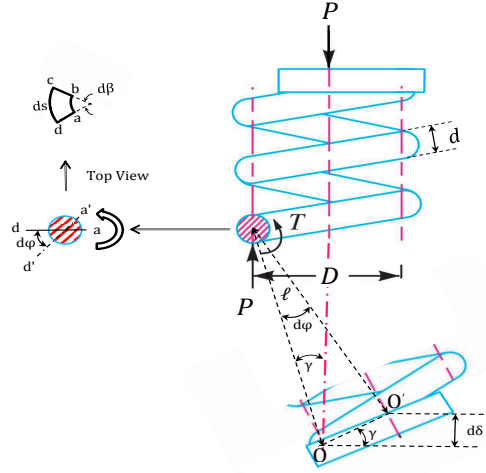, height=2.5 in, width=2.5 in}}
\caption{A Compressive Helical Spring Diagram}
\label{Fig:one}
\end{figure}
We first have
\[GI_n\frac{d\phi}{ds}=\tau ds=P\left(\frac{D}{2}\right)^2\frac{d\beta}{ds}\] which gives $\frac{d\phi}{ds}=\frac{PD^2}{4GI_n}\frac{d\beta}{ds}$. Then, we have \[\frac{d\delta}{ds}=l\sin{\gamma}\frac{d\phi}{ds}
=l\frac{D}{2l}\frac{PD^2}{4GI_n}\frac{d\beta}{ds}\]\[= \frac{(n+3)2^{n+2}}{8}\frac{D^3P}{G\pi d^{n+3}}\frac{d\beta}{ds}\] Therefore, we have the tip deflection formula
\begin{equation}\label{Deflection}
 \delta=\int_0^{2\pi N} \frac{(n+3)2^{n+2}D^3P}{8G\pi d^{n+3}}d\beta= \frac{(n+3)2^{n+2}D^3PN}{4Gd^{n+3}}
\end{equation}
\section{Formulation of the Geometric Programming Problem}\label{sec:3}

The objective function of our optimization problem is the weight of a helical spring under axial load $P$. As for the constraints, we consider only  upper bounds on the  shear stress in the cross-sections of the wire  as described in the previous section and the tip deflection of the spring. First, define the design vector to be
\begin{equation}\label{design_vector}
\mathbf{X} = \left\{\begin{array}{cc}
x_1 \\
x_2
\end{array} \right\} = \left\{\begin{array}{cc}
D \\
d
\end{array} \right\}
\end{equation}
where $D$ is the mean diameter of the coil and $d$ is the diameter of the wire as defined in the nomenclature.
Then, the objective function (mass) of the helical spring can be expressed as
\begin{equation}
f(\mathbf{X}) = \frac{\pi d^2}{4}(\pi D)\rho N
\end{equation}
where $N$ is the number of  turns in the spring and $\rho$ is the density of the spring material. From (\ref{Stress}) with $r=d/2$,  the maximum stress for the power law  helical spring is given by
\[\tau_{max} =  \frac{2(n+3)PD}{\pi d^3}+\frac{4P}{\pi d^2}\]
\begin{equation}
\tau_{max} = \frac{2(n+3)PD}{\pi d^3}K_s
\end{equation}
where $K_s = 1+\frac{2}{(n+3)C}$, $C = \frac{D}{d}$ is called Spring index.
By (3), the expression for the maximum tip deflection of the helical spring is now
\begin{equation}
\delta_{max} = \frac{(n+3)2^{n+3}}{8}\frac{x_1^3 P N}{G x_2^{n+3}}
\end{equation}
where
\begin{equation}
I_n = \frac{2\sqrt{\pi}d^{n+3}}{n+3}\frac{\Gamma(1+\pi/2)}{\Gamma(3/2+\pi/2)} = \frac{\pi d^{n+3}}{(n+3)2^{n+2}}
\end{equation}
We now have the following nonlinear geometric programming problem for the parameters $\rho$, $P$, $\tau_{max}$, $\delta_{max}$, $K$, $n$, $N$:
\begin{equation}\label{prime}
\begin{array}{ll@{}ll}
Minimize & f(x_1,x_2) = \frac{\pi^2 }{4}\rho Nx_2^2x_1 \\ \\
subject \; to & \frac{2(n+3)Px_1}{\pi x_2^3}[1+\frac{2x_2}{(n+3){x_1}}] \leq \tau_{max}\\
           & \hspace{30pt}\frac{(n+3)2^{n+3}}{8}\frac{x_1^3 P N}{Gx_2^{n+3}} \leq \delta_{max}  \\
           & \hspace{88pt} x_2 \leq x_1
\end{array}
\end{equation}
An examination of the KKTC conditions of \eqref{prime} results in no solution. This is because the third constraint must be active i.e. $x_1 = x_2$, which is impossible in the context of the problem.
Therefore, a third constraint, which is  $kx_2 \leq x_1$, where $k>1$ is added to complete the nonlinear program model. This is motivated by using the spring index equation $C=\frac{x_1}{x_2}$ in literature. With the new formulation of the problem our goal is to find optimal values of $k, x_1, x_2$.We have
\begin{equation}\label{new_prime}
\begin{array}{ll@{}ll}
\text{minimize} & f = cx_1x_2^2\\ \\
\text{subject to} &  g_1 = c_{11}x_1x_2^{-3}+c_{12}x_2^{-2}- 1\leq 0 \\
                             & g_2 = c_{21}x_1^3x_2^{-n-3}- 1 \leq 0 \\
                             & g_3 = kx_2 - x_1 \leq 0   \\
                             & x_1, x_2 \geq 0 \\
\end{array}
 \end{equation}
 where $c = \frac{\pi^2\rho N}{4}$, $c_{11} = \frac{2(n+3)P}{\pi \tau_{max}}$, $c_{12} = \frac{4P}{\pi \tau_{max}}$, $c_{21} = \frac{(n+3)2^nPN}{G\delta_{max}}$. Notice that all coefficients are positive.
\section{Computing Optimal Solution from KKTC Conditions}\label{sec:4}
To find optimal values of $k, x_1, x_2$ for the prime problem (\ref{new_prime}), we examine the KKTC conditions as listed below:
\begin{equation*}
\begin{array}{ll@{}ll}
 cx_2^2+c_{11}\lambda_1x_2^{-3}+3c_{21}\lambda_2x_1^2x_2^{-n-3}-\lambda_3 &= 0 \\
 2cx_1x_2-3c_{11}\lambda_1x_1x_2^{-4}-2\lambda_1c_{12}x_2^{-3}- &\\
-\lambda_2(n+3)c_{21}x_1^3x_2^{-n-4}+ k\lambda_3 & = 0 \\
 (c_{11}x_1x_2^{-3}+c_{12}x_2^{-2}-1)\lambda_1  &= 0 \\
 (c_{21}x_1^3x_2^{-n-3}-1)\lambda_2 & = 0 \\
 (kx_2-x_1)\lambda_3 & = 0 \\
 \lambda_i \geq 0 \; \; \text{for} \;\; i = 1,2,3 & \\
 g_i \leq 0 \;\; \text{for} \;\; i=1,2,3 &\\
\end{array}
 \end{equation*}
 If $\lambda_3 = 0$, then from the first equation it follows that $\lambda_1$ or/and $\lambda_2$ are less than $0$. Hence, $\lambda_3$ must be positive.
 \\
 \\
\underline{Case 1: $\lambda_3>0$, $\lambda_1=0$ and $\lambda_2=0$}
\\
\\
 KKTC conditions give
\begin{equation*}
\begin{array}{ll@{}ll}
cx_2^2 - \lambda_3 & = 0 \\
 2cx_1x_2 + k\lambda_3 & = 0 \\
kx_2-x_1 & = 0   \\
\end{array}
 \end{equation*}
 from which $x_2 = 0$ is obtained. However, $x_2$ cannot be zero; therefore, this case results in no solution.
 \\
 \\
 \underline{Case 2: $\lambda_3>0$, $\lambda_1>0$ and $\lambda_2>0$}
 \\
 \\
 From KKTC conditions we have
\begin{equation*}
\begin{array}{ll@{}ll}
 c_{11}x_1x_2^{-3}+c_{12}x_2^{-2}-1 & = 0 \\
 c_{21}x_1^3x_2^{-n-3}-1 & = 0 \\
kx_2-x_1 & = 0   \\
\end{array}
\end{equation*}
This system, in general, has no solution since there are three equations and two unknowns. This case gives no solution, too.
\\
\\
 \underline{Case 3: $\lambda_3>0$, $\lambda_1>0$ and $\lambda_2=0$}
 \\
 \\
 We have the following system of equations
\begin{equation*}
\begin{array}{ll@{}ll}
 cx_2^2+c_{11}\lambda_1x_2^{-3}-\lambda_3 & = 0 \\
 2cx_1x_2-3c_{11}\lambda_1x_1x_2^{-4}-2\lambda_1c_{12}x_2^{-3}+ k\lambda_3 & = 0 \\
 c_{11}x_1x_2^{-3}+c_{12}x_2^{-2}-1 & = 0 \\
 kx_2-x_1 & = 0   \\
\end{array}
\end{equation*}
This gives
\\
\[x_1 = k\sqrt{c_{11}k+c_{12}}, \;\; x_2 = \sqrt{c_{11}k+c_{12}}\]
\[\lambda_1 = \frac{3ckx_2^5}{2c_{11}k+2c_{12}}, \;\; \lambda_3 = \frac{5c_{11}k+2c_{12}}{2c_{11}k+2c_{12}}cx_2^2\]
\\
$\lambda$'s and $x$'s are positive and satisfy non-negativity constraints. In this case, inequality $g_2 \leq 0$ must also hold. That is,
\\
\begin{equation}\label{case3_inq}
{\displaystyle c_{21}x_1^3x_2^{-n-3}-1 \leq 0 \; \Rightarrow \; c_{21}^{2/n}k^{6/n}-c_{11}k-c_{12} \leq 0}
\end{equation}
\\
We must choose such constant $k$ that satisfies the above inequality and minimizes the objective function which becomes
\\
\[f_1 = cx_1x_2^2 = ck(c_{11}k+c_{12})^{3/2}\]
\\
Notice that minimization of the objective function is equivalent to minimization of $k$. Thus, for this case we choose minimum $k>1$ that satisfies~\eqref{case3_inq}.
\\
\\
\underline{Case 4: $\lambda_3>0$, $\lambda_1=0$ and $\lambda_2>0$}
\\
\\
In this case KKTC conditions result in the following system of equations
\begin{equation*}
\begin{array}{ll@{}ll}
 cx_2^2+3c_{21}\lambda_2x_1^2x_2^{-n-3}-\lambda_3 & = 0 \\
 2cx_1x_2-\lambda_2(n+3)c_{21}x_1^3x_2^{-n-4}+ k\lambda_3 & = 0 \\
 c_{21}x_1^3x_2^{-n-3}-1 & = 0 \\
 kx_2-x_1 & = 0   \\
\end{array}
\end{equation*}
from which we obtain
\\
\[x_1 = c_{21}^{1/n}k^{1+3/n}, \;\; x_2 = c_{21}^{1/n}k^{3/n}\]
\[ \lambda_2 = \frac{3cx_2^{n+3}}{c_{21}k^2n}, \;\; \lambda_3 = c\left(1+\frac{9}{n}\right)x_2^2\]
\\
$\lambda$'s and $x$'s satisfy non-negativity constraints. Next, consider inequality $g_1 \leq 0$
\\
\begin{equation}\label{case4_inq}
{\displaystyle c_{11}x_1x_2^{-3}+c_{12}x_2^{-2}-1\leq 0  \Rightarrow c_{21}^{2/n}k^{6/n}-c_{11}k-c_{12} \geq 0}
\end{equation}
\\
Notice that LHS of \eqref{case3_inq} and \eqref{case4_inq} are the same.
\\
\\
In this case, the objective function is
\\
\[f_2 = cx_1x_2^2 = cc_{21}^{3/n}k^{1+9/n}\]
\\
Again, minimization of the objective function is equivalent to minimization it relative to  $k$. Now, we prove that if the following  equation
\\
\begin{equation}\label{k_equation}
{\displaystyle c_{21}^{2/n}k^{6/n}-c_{11}k-c_{12} = 0}
\end{equation}
\\
has a root $k^*>1$ then KKTC conditions provide a desired solution.
Suppose that $k^*$ is a minimal root of \eqref{k_equation} that is greater than $1$. Assume that interval $(1, k^*]$ satisfies inequality \eqref{case3_inq}, then there is no minimal value of $k$ that satisfies the following system of equations\\
\[\left\{
  \begin{array}{rr}
   c_{21}^{2/n}k^{6/n}-c_{11}k-c_{12} \leq 0 & \\
    k>1 & \\
  \end{array}
\right.
\]
Therefore, case 3 does not give any solution. However, case 4 provides a solution because there exists a minimal value of $k$ of the following system of equations
\[\left\{
  \begin{array}{rr}
   c_{21}^{2/n}k^{6/n}-c_{11}k-c_{12} \le 0 & \\
    k>1 & \\
  \end{array}
\right.
\]
and that minimal value of $k$ is $k^*$.
Similarly, if the interval $(1,k^*]$ satisfies inequality \eqref{case4_inq}, then case 3 provides a solution. \\ \\
We now show that solutions in both cases are the same. In other words, $x_1$ and $x_2$ in case 3 are equal to $x_1$ and $x_2$ in case 4, respectively. It is enough to show that $x_2$'s are the same.
\[\sqrt{c_{11}k^*+c_{12}} = c_{21}^{1/n}(k^*)^{3/n} \; \Leftrightarrow \; c_{11}k^*+c_{12} = c_{21}^{2/n}(k^*)^{6/n}\]
\[\Leftrightarrow \; c_{21}^{2/n}(k^*)^{6/n}-c_{11}k^*-c_{12} = 0\]
Since $x_2$'s are equal, $x_1$'s and the values of objective functions in two cases are the same. This means that any solution of $x_1$ and $x_2$ can be chosen from case 3 or 4.
We denote  the function $c_{21}^{2/n}(k^*)^{6/n}-c_{11}k^*-c_{12}$ by $g(k)$.   The algorithm of finding optimal diameters and mass of helical spring can be stated as  the following steps:
\textit{Step 1}: \\
\\
Solve numerically the inequality \\
\[g(k)=c_{21}^{2/n}k^{6/n}-c_{11}k-c_{12} \le 0\]
\\
and choose a value of $k$ in the solution interval of the inequality that is larger than 1 ( See Figure 2),
\\
\textit{Step 2}: \\
\\
Calculate $x_1$, $x_2$ and $f$
\\
\[x_1 = k\sqrt{c_{11}k+c_{12}}, \;\;\; x_2 = \sqrt{c_{11}k+c_{12}}, \;\;\; f = cx_1x_2^2\]

It can be shown ( see also Figure 2) that the equation $g(k)=0$ has one negative root and one positive root that is greater than $1$. This gives optimal design variables for and spring index $C=k$ in the range $(1,k^*]$, where $k^*$  denotes the positive  root.
Notice that we can again optimize the objective function relative to the spring index $C$ by choosing the design variables corresponding to the minimal weight.
\section{Computing the Optimal Solution of the Dual Problem}\label{sec:5}
The algorithm of solving geometric programming problem by constructing its dual is well described in \cite{RAO}. We will follow that algorithm and show that the solution of the dual coincides with the solution of the primal problem.
\\ \\
First, rewrite the problem in the following way
\begin{equation}\label{newprime}
\begin{array}{ll@{}ll}
Minimize & f(x_1,x_2) = cx_1x_2^2 \\ \\
subject \; to & c_{11}x_1x_2^{-3}+c_{12}x_2^{-2} \leq 1\\
           & \hspace{25pt} c_{21}x_1^3x_2^{-(n+3)} \leq 1 \\
           & \hspace{50pt} kx_1^{-1}x_2 \leq 1
\end{array}
\end{equation}
We now form the dual of the primal problem. Remember that the maximum of the dual problem is equal to the minimum of the prime. Hence, we have
\begin{equation}\label{dual}
\begin{array}{ll@{}ll}
Maximize & v = \left( \frac{c}{\lambda_{01}}\lambda_{01}\right)^{\lambda_{01}}\left(\frac{c_{11}}{\lambda_{11}}(\lambda_{11}+\lambda_{12})\right)^{\lambda_{11}} \\
 & \hspace{20pt}\times \left(\frac{c_{12}}{\lambda_{12}}(\lambda_{11}+\lambda_{12})\right)^{\lambda_{12}}\left(\frac{c_{21}}{\lambda_{21}}\lambda_{21}\right)^{\lambda_{21}} \\
 & \hspace{20pt}\times \left( \frac{k}{\lambda_{31}} \lambda_{31}\right)^{\lambda_{31}} \\ \\
subject \; to & \lambda_{01} = 1 \\
           & \lambda_{01}+\lambda_{11}+3\lambda_{21}-\lambda_{31}=0 \\
           & 2\lambda_{01} - 3\lambda_{11} - 2\lambda_{12} - (n+3)\lambda_{21}+\lambda_{31} = 0\\
           & \lambda_{11}+\lambda_{12} \geq 0 \\
           & \lambda_{21} \geq 0 \\
           & \lambda_{31} \geq 0 \\
\end{array}
\end{equation}
From the equations
\\
\\
$\lambda_{01} = 1 $\\
$\lambda_{01}+\lambda_{11}+3\lambda_{21}-\lambda_{31}=0$ \\
$2\lambda_{01} - 3\lambda_{11} - 2\lambda_{12} - (n+3)\lambda_{21}+\lambda_{31} = 0$\\
\\
we obtain
\\
\begin{equation}\label{lambdas}
\begin{array}{ll@{}ll}
\lambda_{21} = \frac{3-2\lambda_{11}-2\lambda_{12}}{n}\\
\lambda_{31} = 1+\frac{9}{n}+(1-\frac{6}{n})\lambda_{11}-\frac{6}{n}\lambda_{12}\\
\end{array}
\end{equation}
\\
Then, the objective function becomes
\begin{equation}\label{dual_obj_function}
\begin{array}{ll@{}ll}
v & \displaystyle = (c)^{\lambda_{01}}\left(\frac{c_{11}}{\lambda_{11}}(\lambda_{11}+\lambda_{12})\right)^{\lambda_{11}}\left(\frac{c_{12}}{\lambda_{12}}(\lambda_{11}+\lambda_{12})\right)^{\lambda_{12}} \\
& \times (c_{21})^{\frac{3-2\lambda_{11}-2\lambda_{12}}{n}}(k)^{1+\frac{9}{n}+(1-\frac{6}{n})\lambda_{11}-\frac{6}{n}\lambda_{12}}
\end{array}
\end{equation}
In order to find maximum of $v$ we differentiate $\ln v$ with repsect to $\lambda_{11}$ and $\lambda_{12}$ and solve the system of two equations.\\
\[\frac{\partial \ln v}{\partial \lambda_{11}} = \ln(c_{11})+\ln(\lambda_{11}+\lambda_{12})-\ln \lambda_{11} - \ln c_{21}^{2/n} + \ln k^{1-6/n}\]
\[\frac{\partial \ln v}{\partial \lambda_{12}} = \ln(c_{12})+\ln(\lambda_{11}+\lambda_{12})-\ln \lambda_{12} - \ln c_{21}^{2/n} - \ln k^{6/n}\]
\\
From
\\
\[\frac{\partial \ln v}{\partial \lambda_{11}} = 0 \;\; \text{and} \;\; \frac{\partial \ln v}{\partial \lambda_{12}} = 0\]
\\
we obtain
\\
\begin{equation}\label{dual_system}
\begin{array}{ll@{}ll}
 \displaystyle \frac{\lambda_{11}+\lambda_{12}}{\lambda_{11}} &= c_{11}^{-1}c_{21}^{2/n}k^{6/n-1}\\ \\
\displaystyle \frac{\lambda_{11}+\lambda_{12}}{\lambda_{12}} &= c_{12}^{-1}c_{21}^{2/n}k^{6/n}\\
\end{array}
\end{equation}
\\
This gives
\\
\[c_{11}^{-1}c_{21}^{2/n}k^{6/n-1} - 1 =\frac{1}{c_{12}^{-1}c_{21}^{2/n}k^{6/n}-1} \Rightarrow \]
\[\Rightarrow c_{21}^{2/n}k^{6/n}-c_{11}k-c_{12} \le 0\]
\\
Notice that the above equation and \eqref{k_equation} are the same.
\\
\\
With \eqref{dual_system} the objective function \eqref{dual_obj_function} becomes
\\
\begin{equation} \label{dual_final_obj_function}
v^* = cc_{21}^{3/n}k^{1+9/n}
\end{equation}
\\
Now, in order to find $x_i$'s we need to solve the following system of equations
\begin{equation}\label{dual_final_system}
\begin{array}{rl@{}rl}
\lambda_{01} & \displaystyle = \frac{cx_1x_2^2}{v^*} \\
\frac{\lambda_{11}}{\lambda_{11}+\lambda_{12}} & = c_{11}x_1x_2^{-3} \\
\frac{\lambda_{12}}{\lambda_{11}+\lambda_{12}} & = c_{12}x_2^{-2} \\
 \frac{\lambda_{21}}{\lambda_{21}} & = c_{21}x_1^3x_2^{-n-3} \\
 \frac{\lambda_{31}}{\lambda_{31}} & = kx_1^{-1}x_2 \\
\end{array}
\end{equation}
\\
With \eqref{dual_system} and \eqref{dual_final_obj_function} the above system becomes
\\
\begin{equation*}
\begin{array}{rl@{}rl}
1 & \displaystyle = \frac{cx_1x_2^2}{cc_{21}^{3/n}k^{1+9/n}} \\
c_{11}^{-1}c_{21}^{2/n}k^{6/n-1} & = c_{11}x_1x_2^{-3} \\
c_{12}^{-1}c_{21}^{2/n}k^{6/n} & = c_{12}x_2^{-2} \\
1 & = c_{21}x_1^3x_2^{-n-3} \\
 1 & = kx_1^{-1}x_2 \\
\end{array}
\end{equation*}
and as a solution we get
\\
\begin{equation} \label{dual_solution}
x_1 = c_{21}^{1/n}k^{1+3/n} \\
x_2 = c_{21}^{1/n}k^{3/n}
\end{equation}
\\
We also need to check whether inequality constraints in \eqref{dual} are satisfied. That is,
\begin{equation*}
\begin{array}{rl@{}rl}
 \lambda_{11}+\lambda_{12} & \geq 0 \\
  \lambda_{21} & \geq 0 \\
  \lambda_{31} & \geq 0 \\
\end{array}
\end{equation*}
 From the first equality in \eqref{lambdas} it follows \\
\begin{equation} \label{sum_lambdas}
\lambda_{21} = \frac{3-2\lambda_{11}-2\lambda_{12}}{n} \geq 0 \Leftrightarrow \lambda_{11}+\lambda_{12} \leq \frac{3}{2}
\end{equation}
\\
Then, from \eqref{dual_system} we obtain \\
\[\lambda_{11}+\lambda_{12} = c_{11}^{-1}c_{21}^{2/n}k^{6/n-1}\lambda_{11} \leq \frac{3}{2}\]
\\
Therefore, \\
\begin{equation}\label{lambda_11}
\lambda_{11} \leq \frac{3}{2}c_{11}c_{21}^{-2/n}k^{1-6/n}
\end{equation}
\\
Combining the second equality in \eqref{lambdas} and \eqref{sum_lambdas} we have
\begin{equation} \label{lambda_31}
\lambda_{31} = 1+\frac{9}{n}+\lambda_{11}-\frac{6}{n}(\lambda_{11}+\lambda_{12}) \geq 1+\lambda_{11} \geq \lambda_{11}
\end{equation}
\\
The values of $\lambda$'s cannot be obtained from \eqref{dual_system} and they do not affect the solution. Hence, we are free to set the values for $\lambda$'s as long as they satisfy the constraints. \\ \\
For $\lambda_{11}$ choose any positive number that satisifies \eqref{lambda_11}. Then, \eqref{dual_system} gives \\
\[\lambda_{12} = (c_{11}^{-1}c_{21}^{2/n}k^{6/n-1}-1)\lambda_{11}\]
\\
From \eqref{sum_lambdas} and \eqref{lambda_31} it follows that $\lambda_{21}$ and $\lambda_{31}$ are positive.
The solution \eqref{dual_solution} and the objective function \eqref{dual_final_obj_function} coincide with the solution in case 3 of the primal problem. Remember that the solutions obtained from case 3 and 4 are the same. The coefficient $k$ is again found by solving the equation $c_{21}^{2/n}k^{6/n}-c_{11}k-c_{12} = 0$. From \eqref{dual_final_obj_function} one should not think that because we maximize $v$ we need to maximize $k$. In the dual and primal problems $k$ is a coefficient and the value of which are free to set. From the solution of the dual we know that $k$ must satisfy \eqref{k_equation}. However, in order to choose the right $k$ one needs to look at the primal problem. In the context of the primal problem $k$ should be minimized i.e. must be the minimum positive  root of \eqref{k_equation} which is greater than 1.
\section{Numerical Experiments}\label{sec:6}
\textbf{Problem:} Formulate the problem of minimum weight design if a helical spring under axial load that is made of stainless steel. Consider constraints on the shear stress and the deflection of the spring. Number of active turns $N = 10$, the density $\rho = 7700 kg/m^3$, $n = 0.1$ for $\nu = 0.275$, $K = 960 MPa$, for the axial load, take $P = 10.0 N$. The maximum deflection of the spring $\delta_{max} = 0.03 m$ and maximum shear stress is $\tau_{max} = 200 MPa$.
First, we numerically solve for minimum root of \eqref{k_equation} and then obtain the optimal solution $x_1,x_2$ from the solution obtained in the previous two sections.
After computing on \texttt{\large{Mathematica}}, we have
$c=189989.8847, c_{11}=9.8676\cdot 10^{-8},c_{12}=6.3662\cdot 10^{-8},c_{21}=1.4709\cdot 10^{-5}$.
The solution of the inequality\eqref{k_equation} is the interval$\quad(-0.645162,33.0756].$
\begin{figure}[ht]
\centerline{\epsfig{figure=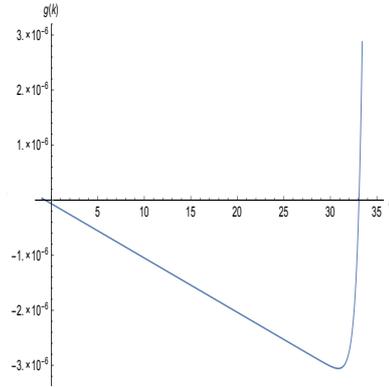, height=2.0 in, width=2.0 in}}
\caption{A Compressive Helical Spring Diagram}
\label{fig:two}
\end{figure}
We now form the prime problem with numerical coefficients for this example:
\begin{equation} \label{eq:er1.7}
	\begin{array}{rl}
 \textnormal{Minimize} & \displaystyle f(x_1, x_2) = 189989.8847 x_1 x_2^2 \\
\textnormal{Subject to} & \displaystyle  \frac{9.8676\cdot 10^{-8} x_1 + 6.3662\cdot 10^{-8} x_2}{x_2^3} \leq 1\\
		& \displaystyle  \frac{1.4709\cdot 10^{-5} x_1^3}{x_2^{3.1}} \leq 1.\\
		& \displaystyle kx_2\le x_1\\
		& \displaystyle x_2, x_1\ge 0
	\end{array}
\end{equation}
where the value of $k$ can be chosen from the interval $(1,33.0756]$.  By solving this prime problem directly and by the semi-analytic solution provided in the steps in section 5, the numerical solutions match with little error. For example, by taking $k=10$, we have the solution $x_1 =0.010249 m$, $x_2 = 0.0010249 m$ and the objective function value is $f = 0.00306809 kg$. This numerical experiment indicates that the solution obtained using the formula in Section 5 agrees with the solutions computed by the \texttt{\large{Mathematica}}.
\section{Conclusion} \label{sec:7}
A nonlinear geometric programing problem is formulated for compression helical springs made of power law materials. Both the prime and the dual problems are shown to have the same solution by examining the KKTC conditions. A semi-analytic solution is derived which provides solutions of the helical spring for a range of spring index.  A numerical example is also provided to illustrate and validate the semi-analytic solution with the  solution computed by solving the prime problem numerically using \texttt{\large{Mathematica}}.


%


\bibliographystyle{spmpsci}      


\end{document}